\documentclass{article}
\input{tcilatex}
\begin{document}

\bigskip 

\bigskip 

\begin{center}
\bigskip {\large A note regarding some generalization of a Schur's theorem}

Igor Ya. Subbotin, PhD,

Department of Mathematics and Natural Sciences, 

College of Letters and Sciences, 

National University, 5245 Pacific Concourse Drive, 

Los Angeles, CA 90045-6904

isubboti@nu.edu

\bigskip 
\end{center}

The article ON THE UPPER CENTRAL SERIES OF INFINITE GROUPS by M. DE FALCO,
F. DE GIOVANNI, C. MUSELLA, AND Y. P. SYSAK, PROCEEDINGS OF THE AMERICAN
MATHEMATICAL SOCIETY, Volume 139, Number 2, February 2011, 385--389 consists
of a quite long proof of the following theorem. 

\textbf{Theorem 3.5} (M. De Falco - F. de Giovanni - C. Musella - Y.P.
Sysak). \ \textit{Let }$G$\textit{\ be a group. The factor group }$G/\bar{Z}%
(G)$\textit{\ is finite if and only if there exists a finite normal subgroup 
}$N$\textit{\ of }$G$\textit{\ such that }$G/N$\textit{\ is
hypercentral.\bigskip }

We can offer the following simple and very brief proof of this result above.

\bigskip 

\textbf{Proof}. Put $Z=\bar{Z}(G)$. Observe at ones that if the length $zl(G)
$ of the upper central series of $G$ is finite, the result follows from the
Baer's theorem. So we will suppose that the length of the upper central
series is infinite. Let $F$ be a finitely generated subgroup with the
property $G=ZF$. Since $G$ is hypercentral-by-finite, every finitely
generated subgroup of $G$ is nilpotent-by-finite, in particular $F$ is
nilpotent-by-finite. We note that $G/Z$ is non-nilpotent, so that $F$ is
non-nilpotent. Let $t=zl(F)$ and $C$ be the upper hypercenter of $F$. If we
suppose that $C\neq C\cap Z$, then $C/(C\cap Z)$ is not trivial. Then $%
CZ/Z\neq \left\langle 1\right\rangle $, which means that the upper
hypercenter of $G/Z$ is not trivial. This contradiction shows that $C\leq Z$%
. By the Baer's theorem, $\gamma _{t+1}(F)$ is finite. Let $L$ be the
nilpotent residual of $F$. Then $L$ is finite. Let $d$ be the nilpotence
class of $F/L$. Choose the local family $\mathcal{M}$ of finitely generated
subgroups of $G$ each of which contains $F$. Moreover, since $zl(G)$ is
infinite, we can choose $\mathcal{M}$ such that $zl(K)\geq d$ for every $%
K\in \mathcal{M}$. Let $K\in \mathcal{M}$ and $C_{K}$ be the upper
hypercenter of $K$. From $G=ZK$ it follows that $C_{K}\leq Z$. We have $%
K=F(K\cap Z)=FC_{K}$. Further, $k=zl(K)\geq d$. The equation $K=C_{K}F$
implies that $\gamma _{k+1}(K)=\gamma _{k+1}(F)$ (see, for example, [N.S.
HEKSTER. ON THE STRUCTURE OF $n$-ISOCLINISM CLASSES OF GROUPS. Journal of
Pure and Applied Algebra 40 (1986), 63-85], Lemma 2.4). The choice of $k$
yields that $L=\gamma _{k+1}(F)$. It follows that $L$ is normal in $K$.
Since $G=\bigcup\limits_{K\in \mathcal{M}}K$, $L$ is normal in $G$.
Furthermore, $K/L$ is nilpotent. Hence $G/L$ has a local family of nilpotent
subgroup. It follows that $G/L$ is locally nilpotent. Being
hypercentral-by-finite and locally nilpotent, $G/L$ is hypercentral.\bigskip 

\end{document}